\documentclass[12pt]{article}
\usepackage{latexsym}
\usepackage{amsmath,amssymb}
\usepackage{bm}
\usepackage{graphicx,xcolor}
\usepackage{amssymb}
\usepackage{cite}
\usepackage{amsfonts, amscd, amsthm}
\usepackage[all]{xy}
\usepackage{here}
\usepackage{ulem}
\usepackage{mathabx}
\usepackage[titletoc,title]{appendix}
\usepackage{mathbbol}

\newcommand{\sn}{{\rm sn}}

\renewcommand{\theequation}{\thesection.\arabic{equation}}

\newcommand{\wwp}{\widehat{\wp}}
\newcommand{\Q}{\widehat{Q}}
\newcommand{\lam}{\lambda}

\makeatletter
\@addtoreset{equation}{section}
\renewcommand{\theequation}{\thesection.\@arabic\c@equation}
\makeatother

\makeatletter
\renewcommand\appendix{\par
  \setcounter{section}{0}%
  \setcounter{subsection}{0}%
  \gdef\thesection{Appendix \@Alph\c@section }
  \renewcommand{\theequation}
  {\Alph{section}.\arabic{equation}}
}
\makeatother

\makeatletter
\newcounter{subeqncnt}
\def\thesubeqncnt{\alph{subeqncnt}}
\def\subequations{\begingroup%
\stepcounter{equation}\edef\@tempa{\theequation}%
\let\c@equation\c@subeqncnt\c@subeqncnt\z@
\edef\theequation{\@tempa\noexpand\thesubeqncnt}}

\makeatother

\allowdisplaybreaks

\begin{document}

\titlepage

\title{Generalized mKdV Equation and Genus Two Jacobi Type 
Hyperelliptic Differential Equation} 
\author{Masahito Hayashi\thanks{masahito.hayashi@oit.ac.jp}\\
Osaka Institute of Technology, Osaka 535-8585, Japan\\
Kazuyasu Shigemoto\thanks{shigemot@tezukayama-u.ac.jp} \\
Tezukayama University, Nara 631-8501, Japan\\
Takuya Tsukioka\thanks{tsukioka@bukkyo-u.ac.jp}\\
Bukkyo University, Kyoto 603-8301, Japan\\
}
\date{\empty}


\maketitle

\abstract{%
We generalized the mKdV equation such that the static equations included the 
$\sn$ differential equation. 
As a result, a good correspondence was obtained between the KdV and mKdV equations.
For the general genus two hyperelliptic curves, we obtained differential equations 
for the Weierstrass-type 
and Jacobi-type hyperelliptic functions.
Considering the special case of $\lambda_6=0, \lambda_0=0$,
the Weierstrass-type and Jacobi-type hyperelliptic functions are different solutions 
to the same hyperelliptic differential equations. These solutions are then connected via a  
special ${\rm Sp(4, {\bf R})}$ Lie group transformation.
 
}


\section{Introduction} 
\setcounter{equation}{0}

The soliton system can be solved exactly despite being nonlinear because of 
the Lie group structure. 
It is well known that the KdV equation has an ${\rm Sp(2,{\bf R})/Z_2}\cong {\rm SO(2,1)}$
Lie group structure~\cite{Crampin,Sasaki}. 
We have demonstrated that $\wp$ function, which is the special static solution of the KdV 
equation,  has the ${\rm Sp(2,{\bf R})/Z_2}\cong {\rm SO(2,1)}$
Lie group structure~\cite{H1}, and the genus two hyperelliptic $\wp$ functions
have the ${\rm Sp(4,{\bf R})/Z_2}\cong {\rm SO(3,2)}$ Lie group structure~\cite{H1,H2,H3}.
In the hyperelliptic curves, not in general Riemann surfaces, we expect that 
multi-periodic functions, that is, hyperelliptic functions, have a Lie group structure.
 
In this study, we first generalized the mKdV equation such that  the static differential equation 
included the $\sn$ differential equation. This gives a good correspondence between the KdV 
and mKdV equations. 
For the general genus two hyperelliptic curves, we obtain the conventional differential equations 
for the Weierstrass-type and Jacobi-type hyperelliptic functions.
If we consider the special case of $\lambda_6=0, \lambda_0=0$, differential equations for 
the Weierstrass-type and Jacobi-type hyperelliptic functions are different solutions of 
the same hyperelliptic differential equations.
We confirmed that the two hyperelliptic functions were linked by a special 
half-period ${\rm Sp(4, {\bf R})}$ Lie group transformation.

\section{Generalized mKdV equation} 
\setcounter{equation}{0}
%
\subsection{Integrability condition for the generalized mKdV equation}
The KdV equation is expressed as follows
\begin{equation}
u_t(x,t)+u_{xxx}(x,t)-6 u(x,t) u_x(x,t)=0,
\label{2e1}
\end{equation}
and mKdV equation is expressed as follows
\begin{equation}
v_t(x,t)+v_{xxx}(x,t)-6 v^2(x,t) v_x(x,t)=0,
\label{2e2}
\end{equation}
where we denote
$\dfrac{\partial u(x,t)}{\partial t}=u_t(x,t)$, 
$\dfrac{\partial u(x,t)}{\partial x}=u_x(x,t)$, 
$\dfrac{\partial^3 u(x,t)}{\partial x^3}=u_{xxx}(x,t)\rule[-3.5mm]{0mm}{9mm}$ etc. 
Eq.(\ref{2e2}) is generalized such that the static equation includes a differential 
equation for the $\sn$ function as follows 
\begin{eqnarray}
v_t(x,t)+v_{xxx}(x,t)-6 v^2(x,t) v_x(x,t)+a v_x(x,t)=0.
\label{2e3}
\end{eqnarray}
We call Eq.(\ref{2e3}) as the generalized mKdV equation.
We demonstrate that this generalized mKdV equation is an integrable system. 
To achieve this, we adopt the following AKNS formalism~\cite{AKNS,Wadati} 
\begin{eqnarray}
&&\frac{\partial}{\partial x} 
\left( \begin{array}{c} \psi_1 \\ \psi_2 
\end{array} \right)
= \left(\begin{array}{cc} \eta & v(x,t) \\ -v(x,t) & -\eta
\end{array}\right)
 \left(\begin{array}{c} \psi_1 \\ \psi_2
\end{array}\right)
=L\left( \begin{array}{c} \psi_1 \\ \psi_2 
\end{array} \right) ,
\label{2e4} \\
&&\frac{\partial}{\partial t} 
\left( \begin{array}{c} \psi_1 \\ \psi_2 
\end{array} \right)
= \left(\begin{array}{cc} A & B \\ C & -A
\end{array}\right)
 \left(\begin{array}{c} \psi_1 \\ \psi_2
\end{array}\right)
=M\left( \begin{array}{c} \psi_1 \\ \psi_2 
\end{array} \right), 
\label{2e5}\\
&&A=-2 \eta v^2(x,t) -4 \eta^3 b, 
\label{2e6}\\
&&B=-v_{xx}(x,t)-2 \eta v_x(x,t) -2 v^3(x,t)-4 \eta^2b v(x,t) ,
\label{2e7}\\
&& C=v_{xx}(x,t)-2 \eta v_x(x,t)+2 v^3(x,t)+4 \eta^2b v(x,t) .
\label{2e8}
\end{eqnarray}
$\eta$ is the spectral parameter and $b$ is a constant 
parameter, which is determined later.  The matrix $M$ in Eq.(\ref{2e5}) has 
${\rm GL(2, {\bf R})}$ structure.
The integrability condition for Eq.(\ref{2e4}) and (\ref{2e5}) are given by
\begin{eqnarray}
&&\left[ \frac{\partial}{\partial x} -L,  \frac{\partial}{\partial t} -M\right]
=\frac{\partial L}{\partial t}-\frac{\partial M}{\partial x}
+\left[ L, M\right]
=\left(\begin{array}{cc} 0 & D \\ -D & 0
\end{array}\right)=0, 
\label{2e9}\\
&&D=v_t(x,t)+v_{xxx}(x,t)+6 v^2(x,t) v_{x}(x,t)+4 \eta^2 (b-1)v_{x}(x,t),
\label{2e10}
\end{eqnarray}
and $b$ is set to be $a=4 \eta^2 (b-1)$. If we replace $v(x,t) \rightarrow i v(x,t)$, the above integrability
condition $D=0$ yields the generalized mKdV equation.

\subsection{Various transformations between the static KdV and the static generalized mKdV equations}
For simplicity, we will denote $v(x,0)$ as $v(x)$. This type of abbreviation is also
 used in other cases.
\subsubsection{Miura transformation}

Using the static Miura transformation $u(x)=v^2(x)+v_x(x)-a/6$, we obtain 
\begin{align}
&u_{xxx}(x)-6 u(x) u_x(x)
=\left(\frac{d}{d x} +2v(x)\right)
\left(v_{xxx}(x)-6 v^2(x) v_x(x)+a v_x(x)\right) .
\label{2e11}
\end{align}
Therefore, if
\begin{equation}
  v_{xxx}(x)-6 v^2(x) v_x(x)+a v_x(x)=0,
\label{2e12}
\end{equation}
we obtain the static KdV equation
\begin{equation}
  u_{xxx}(x)-6 u(x) u_x(x)=0.
\label{2e13}
\end{equation}
The time-dependent KdV and generalized mKdV equations are linked by the time-dependent Miura transformation $u(x,t)=v^2(x,t)+v_x(x,t)-a/6$~\cite{Miura}.  

\subsubsection{Square transformation}

Using the square transformation $u(x)=2 v^2(x)-2a/3$, we obtain 
\begin{align}
&u_{xxx}(x)-6 u(x) u_x(x)
=\left(4 v(x)\frac{d}{d x} +12v_x(x)\right)
\left(v_{xx}(x)-2 v^3(x)+a v(x)\right) .
\label{2e14}
\end{align}
Therefore, if
\begin{equation}
  v_{xx}(x)-2 v^3(x)+a v(x)=0,
\label{2e15}
\end{equation}
we obtain the static KdV equation
\begin{equation}
  u_{xxx}(x)-6 u(x) u_x(x)=0.
\label{2e16}
\end{equation}
Note that
\begin{equation}
  \frac{d}{dx}\left(v_{xx}(x)-2 v^3(x)+a v(x)\right)
  =
  v_{xxx}(x)-6 v^2(x) v_x(x)+a v_x(x).
  \label{2e17}
\end{equation}
Eq.\eqref{2e15} means this $v(x)$ satisfies the static generalized mKdV equation.

\subsubsection{Inverse square transformation}

Using the inverse square transformation $u(x)=1/v^2(x)-2a/3$, we obtain 
\begin{align}
u_{xxx}&(x)-6 u(x) u_x(x)
\nonumber\\
=&\left( -\frac{1}{v^3(x) v_x(x)}  \frac{d^2}{d x^2}
+\left(\frac{v_{xx}(x)}{v^3(x) {v_x^2(x)}}+\frac{9}{v^4(x)}\right)
\frac{d}{d x}
-\frac{24 v_x(x)}{v^5(x)} \right)
\nonumber\\
& \times
\left(v_{x}^2(x)-v^4(x)+a v^2(x)-\frac{1}{2} \right) .
\label{2e18}
\end{align}
Therefore, if
\begin{equation}
  v_{x}^2(x)-v^4(x) +a v^2(x)-\frac12 =0,
\label{2e19}
\end{equation}
we obtain  the static KdV equation
$$u_{xxx}(x)-6 u(x) u_x(x)=0.$$
Note that
\begin{equation}
  \frac{d}{dx}\left(v_{x}^2(x)-v^4(x) +a v^2(x)-\frac12\right)
  =
  2v_x(x)
  \left(v_{xx}(x)-2 v^3(x)+a v(x)\right).
\label{2e20}
\end{equation}
In this case, $v(x)$ not only satisfies the static generalized mKdV equation,  
but also becomes the Jacobi $\sn$ function.
This transformation is well-known to link
the Weierstrass $\wp$ function and Jacobi 
$\sn$ function in the form  
$$\wp\left(\frac{u}{\sqrt{e_1-e_3}}\right)=e_3+\frac{(e_1-e_3)}{\sn^2(u)}.$$

\subsubsection{Inverse power transformation}

Using the inverse power transformation $u(x)=1/v(x)-a/6$, we obtain 
\begin{align}
u_{xxx}&(x)-6 u(x) u_x(x)
\nonumber\\
=&\left( -\frac{1}{2 v^2(x) v_x(x)}  \frac{d^2}{d x^2}
+\left(\frac{v_{xx}(x)}{2 v^2(x) {v_x(x)}^2}+\frac{3}{v^3(x)}\right)
\frac{d}{d x}
-\frac{6v_x(x)}{v^4(x)} \right)
\nonumber\\
& \times
\left(v_{x}^2(x)-v^4(x) +a v^2(x)-2 v(x)\right) .
\label{2e21}
\end{align}
Therefore, if
\begin{equation}
  v_{x}^2(x)-v^4(x) +a v^2(x)-2 v(x)=0,
\label{2e22}
\end{equation}
we obtain the static KdV equation
$$u_{xxx}(x)-6 u(x) u_x(x)=0.$$
Note that
$$
\frac{d}{dx}\left(v_{x}^2(x)-v^4(x) +a v^2(x)-2 v(x)\right)
=
2v_x(x)
\left(v_{xx}(x)-2 v^3(x)+a v(x)-1\right),
$$
and
$$
\frac{d}{dx}\left(v_{xx}(x)-2 v^3(x)+a v(x)-1\right)
= 
v_{xxx}(x)-6v^2(x)v_x(x)+av_x(x).
$$
In this case, $v(x)$ is a solution of the static generalized mKdV equation, but not 
that of the Jacobi $\sn$ function.
\subsubsection{A relationship between transformations}
Comparing the square transformation and the inverse square 
transformation, we obtain
\begin{equation} 
u(x)=2{v_1^2(x)}-\frac{2a}{3}=\frac{1}{{v_2^2(x)}}-\frac{2a}{3},
\label{2e23}
\end{equation}
where $v_1(x)$ and $v_2(x)$ satisfy the same $\sn$-type differential equation. 
In order that Eq.(\ref{2e22}) becomes the standard $\sn$ differential equation, we set
$k=\sqrt{2}$, $a=(1+k^2)/2=3/2$, and 
obtain
\begin{equation}
v_1(x)=\sn(x/\sqrt{2}) ~~\mathrm{and}~~
 v_2(x)=\widehat{\sn}(x/\sqrt{2})
=\frac{1}{\sqrt{2}\ \sn(x/\sqrt{2})} =\frac{1}{\sqrt{2}\ v_1(x)}.
\label{2e24}
\end{equation}
More generally, using 
$\sn_z^2(z)=\left(1-\sn^2(z)\right)\left(1-k^2 \sn^2(z)\right)$,
the differential equation for $\sn(z)$ becomes
\begin{eqnarray}
\sn_{zz}(z)=-\left(1+k^2\right) \sn(z)+2 k^2 \sn^3(z) .
\label{2e25}
\end{eqnarray}
On the other hand, $\displaystyle{\widehat{\sn}(z)=\frac{1}{k\,\sn(z)}}$ satisfies the same differential equation~\cite{H1}
\begin{eqnarray}
\widehat{\sn}_{zz}(z)=-\left(1+k^2\right) \widehat{\sn}(z)+2 k^2 \widehat{\sn}^3(z).
\label{2e26}
\end{eqnarray}
The reason for this relation comes from the half-period relation 
$\displaystyle{\sn(z+3 i K')=\frac{1}{k\ \sn(z)}}$.
\section{Differential equations for genus two Jacobi type hyperelliptic functions} 
\setcounter{equation}{0}

We start from the genus two hyperelliptic curves of the form
\begin{eqnarray}
{y_i}^2=\lambda_6 {x_i}^6+\lambda_5 {x_i}^5+\lambda_4 {x_i}^4+\lambda_3 {x_i}^3 
+\lambda_2 {x_i}^2+\lambda_1 {x_i}+\lambda_0, \quad (i=1, 2). 
\label{3e1} 
\end{eqnarray}
The Jacobi inversion relation is given by 
\begin{eqnarray}
du_1=\frac{dx_1}{y_1}+\frac{dx_2}{y_2}, \quad 
du_2=\frac{x_1 dx_1}{y_1}+\frac{x_2 dx_2}{y_2},  
\label{3e2}
\end{eqnarray}
which gives  
\begin{eqnarray}
\frac{\partial x_1}{\partial u_2}=\frac{y_1}{x_1-x_2}, \quad
\frac{\partial x_2}{\partial u_2}=-\frac{y_2}{x_1-x_2}, \quad
\frac{\partial x_1}{\partial u_1}=-\frac{x_2 y_1}{x_1-x_2}, \quad
\frac{\partial x_2}{\partial u_1}=\frac{x_1 y_2}{x_1-x_2} .
\label{3e3}
\end{eqnarray}
We define
\begin{align}
F(x_1,x_2; \left\{\lambda_j\right\})
=&2 \lambda_6 x_1^3 x_2^3+\lambda_5 x_1^2 x_2^2 (x_1+x_2) 
+2 \lambda_4 x_1^2 x_2^2+\lambda_3 x_1 x_2 (x_1+x_2)
\nonumber\\
&+2 \lambda_2 x_1 x_2+\lambda_1 (x_1+x_2)+2 \lambda_0.
\label{3e4}
\end{align}

Then the genus two Weierstrass type hyperelliptic functions are defined 
as follows~\cite{Baker1} 
\begin{eqnarray}
&&\Re_{22}(u_1, u_2)=\frac{\lambda_5}{4}
\left( x_1+x_2+\frac{2 \lambda_6}{\lambda_5}(x_1^2+x_1 x_2+x_2^2)\right), 
\label{3e5}\\
&&\Re_{21}(u_1, u_2)=\frac{\lambda_5}{4}
\left( -x_1x_2 -\frac{2 \lambda_6}{\lambda_5}x_1 x_2 (x_1+x_2)\right), 
\label{3e6}\\
&&\Re_{11}(u_1, u_2)=\frac{F(x_1,x_2;\left\{\lambda_j\right\})-2 y_1 y_2}{4 (x_1-x_2)^2}
+ \frac{\lambda_6}{2} x_1^2 x_2^2 .
\label{3e7}
\end{eqnarray}
where we denote $\displaystyle{\frac{\partial^2 \Re}{\partial u_i \partial u_j}\rule[-4mm]{0mm}{9mm}= 
\Re_{i j}}$, $\displaystyle{\frac{\partial^4 \Re}{\partial u_i \partial u_j \partial u_k \partial u_l}= 
\Re_{i j k l}}$, $(i, j, k, l=1, 2)$ etc.
These hyperelliptic functions satisfy the complete integrability of the form
\begin{eqnarray}
\frac{ \partial \Re_{22}}{ \partial u_1}
=\frac{ \partial \Re_{21}}{ \partial u_2}, \quad
\frac{ \partial \Re_{21}}{ \partial u_1}
= \frac{ \partial \Re_{11}}{ \partial u_2}.
\label{3e8}
\end{eqnarray}
Though $\Re_{22}$ and $\Re_{21}$ are symmetric functions,
but they are not the fundamental symmetric functions. 
Therefore, it is not always possible to express any symmetric function 
in a closed form using $\Re_{22}$ and $\Re_{21}$.
Until now, closed differential equations for the genus two hyperelliptic
 functions are known only for the $\lambda_6=0$ case~\cite{H4}.

\subsection{Genus two Weierstrass type hyperelliptic differential equations}
Even for the general genus two hyperelliptic curves with $\lambda_6 \ne 0$, 
we consider the genus two hyperelliptic functions  
$\wp_{22}$, $\wp_{21}$ and $Q$ of the form 
\begin{align}
\wp_{22}(u_1, u_2)=\frac{\lambda_5}{4} ( x_1+x_2), 
\ \wp_{21}(u_1, u_2)=-\frac{\lambda_5}{4} x_1x_2, 
\ Q(u_1, u_2)=\frac{F(x_1,x_2;\left\{\lambda_j\right\})-2 y_1 y_2}{4 (x_1-x_2)^2}.
\label{3e9}
\end{align}
We refer to this as Weierstrass-type hyperelliptic functions. 
Here $\wp_{22}$ and $\wp_{21}$ are fundamental symmetric function. Then any symmetric 
quantity can be expressed by $\wp_{22}$ and $\wp_{21}$.
Thus, we can obtain hyperelliptic differential equations in closed forms.
Though we obtain the relation
\begin{eqnarray}
\frac{ \partial \wp_{22}}{ \partial u_1}
=\frac{ \partial \wp_{21}}{ \partial u_2}, \quad
\frac{ \partial \wp_{21}}{ \partial u_1}
\ne \frac{ \partial Q}{ \partial u_2} ,
\label{3e10}
\end{eqnarray}
the condition 
$\displaystyle{\frac{ \partial \wp_{22}}{ \partial u_1}
=\frac{ \partial \wp_{21}}{ \partial u_2}}$ 
alone is sufficient to guarantee the integrability condition, that is, the existence of $\sigma(u_1,u_2)$ function as a potential. 
If $\wp_{22}$ and $\wp_{21}$ are given as functions of $\sigma$ function, $Q$ can in principle 
be given as a function of $\sigma$ function via the generalized Kummer surface relation. 
We will provide a generalized Kummer surface relation later.

In order to obtain the differential equations for genus two Weierstrass type hyperelliptic
functions, we first calculate 
$\wp_{2222}$, $\wp_{2221}$, $\wp_{2211}$, $\wp_{2111}$, $Q_{11}$, $Q_{21}$ and $Q_{22}$. 
Because these quantities are symmetric functions of $x_1$ and $x_2$, these are 
expressed by $\wp_{22}$ and $\wp_{21}$. If these quantities contain symmetric function 
of the type $y_1 y_2$, the quantity $y_1 y_2$ is expressed  by $Q$ and/or $Q^2$ by using 
$\displaystyle{y_1 y_2=-2(x_1-x_2)^2 Q+ \frac{1}{2} F(x_1,x_2)}$. 

Then we obtain the differential equations of the genus two Weierstrass type hyperelliptic
functions of the form 
\begin{align}
&1)\ \wp_{2222}- \frac{32\lambda_6}{\lambda_5^2} \wp_{22}^3
-\frac{12 \lambda_6}{\lambda_5} \wp_{22}\wp_{21}
-6 \wp_{22}^2-\lambda_4 \wp_{22}- \lambda_5 \wp_{21}
-\frac{\lambda_3 \lambda_5}{8} 
 =0,
\label{3e11}\\
&2)\ \wp_{2221}
- \frac{32\lambda_6}{\lambda_5^2} \wp_{22}^2 \wp_{21}
-\frac{4 \lambda_6}{\lambda_5} \wp_{21}^2
-6 \wp_{22}\wp_{21}-\lambda_4 \wp_{21}+\frac{\lambda_5}{2} Q
=0,
\label{3e12}\\
&3)\ \wp_{2211}
- \frac{32\lambda_6}{\lambda_5^2} \wp_{22} \wp_{21}^2
-2 \wp_{22} Q-4 \wp_{21}^2
-\frac{\lambda_3}{2} \wp_{21}=0,
\label{3e13}\\
&4)\ \wp_{2111}
- \frac{32\lambda_6}{\lambda_5^2} \wp_{21}^3
-6 \wp_{21} Q
+\frac{\lambda_1}{2} \wp_{22}-\lambda_2 \wp_{21}+\frac{\lambda_0 \lambda_5}{4}=0,
\label{3e14}\\
&5)\ Q_{11}
- \frac{32\lambda_6}{\lambda_5^2} \wp_{21}^2 Q
+ \frac{16 \lambda_0 \lambda_6}{\lambda_5^2} \wp_{22}^2
- \frac{8 \lambda_1 \lambda_6}{\lambda_5^2} \wp_{22}\wp_{21}
-6 Q^2 +3 \lambda_0 \wp_{22}
\nonumber\\
&-\left(\lambda_1-\frac{2 \lambda_0 \lambda_6}{\lambda_5}\right) \wp_{21}-\lambda_2 Q
+\frac{\lambda_0 \lambda_4}{2}-\frac{\lambda_1 \lambda_3}{8}=0,
\label{3e15}\\
&6)\ Q_{21}
 -\frac{32 \lambda_6}{\lambda_5^2} \wp_{22} \wp_{21} Q
+\frac{4 \lambda_1 \lambda_6}{\lambda_5^2} \wp_{22}^2
-\frac{8 \lambda_2 \lambda_6}{\lambda_5^2} \wp_{22} \wp_{21}
+\frac{4 \lambda_3 \lambda_6}{\lambda_5^2} \wp_{21}^2
\nonumber\\
&-6 \wp_{21} Q 
+\left(\frac{2 \lambda_0 \lambda_6}{\lambda_5} +\frac{\lambda_1}{2}\right)\wp_{22}
-\lambda_2 \wp_{21}+\frac{\lambda_0 \lambda_5}{4}=0,
\label{3e16}\\
&7)\ Q_{22}
-\frac{32 \lambda_6}{\lambda_5^2} \wp_{22} ^2 Q
+ \left( \frac{16 \lambda_4 \lambda_6}{\lambda_5^2} -4\right) \wp_{21}^2
-\frac{8 \lambda_3 \lambda_6}{\lambda_5^2} \wp_{22} \wp_{21}
-2 \wp_{22} Q
\nonumber\\
&
-\frac{20 \lambda_6}{\lambda_5} \wp_{21} Q
+\frac{\lambda_1 \lambda_6}{\lambda_5} \wp_{22}
+\left(-\frac{2 \lambda_2 \lambda_6}{\lambda_5}- \frac{\lambda_3}{2}\right)\wp_{21}
+\frac{\lambda_0 \lambda_6}{2}=0.
\label{3e17}
\end{align}
From Eq.(\ref{3e12}), $Q$ naturally appears as a combination of the terms in Eq.(\ref{3e9}).
The characteristic feature of $\lambda_6 \ne 0$ hyperelliptic differential equations is the existence
 of the cubic term, that is, terms 
$\displaystyle{-\frac{32 \lambda_6}{\lambda_5^2} \wp_{22}^3,
-\frac{32 \lambda_6}{\lambda_5^2} \wp_{22}^2 \wp_{21}, 
 \cdots, -\frac{32 \lambda_6}{\lambda_5^2} \wp_{22}^2 Q}$,  
 in Eq.({\ref{3e11}) $\sim$ Eq.(\ref{3e17}).\\

Here we present a generalized Kummer surface relation.
For $\wp_{22}$, $\wp_{21}$ and $Q$, only two are 
independent. Then, there exists one relation, called the generalized Kummer surface relation.
Using the relation 
$\displaystyle{y_1 y_2=-2(x_1-x_2)^2 Q+ \frac{1}{2} F(x_1,x_2)}$, we obtain 
$$\displaystyle{y_1^2 y_2^2=4(x_1-x_2)^4 Q^2-2(x_1-x_2)^2 Q F(x_1,x_2)+\frac{1}{4} F(x_1,x_2)^2}.$$ 
Using Eq.(\ref{3e1}), this constraint is expressed by $\wp_{22}$, $\wp_{21}$ and $Q$.
The explicit expression is as follows.
We denote $X=\wp_{22}$, $Y=\wp_{21}$, $Z=Q$, and obtain the following 
generalized Kummer surface relation of the form 
\begin{align}
&K_2=K_1
+\frac{4\lambda_6}{\lambda_5^2}
\Big(-\lambda_5^2 Y^2 -8 \lambda_0 \lambda_5 X^2 Y
+4 \lambda_1 \lambda_5 X Y^2
\nonumber\\
&+16(-\lambda_0 X^4+\lambda_1 X^3 Y-\lambda_2 X^2 Y^2
+\lambda_3 X Y^3-\lambda_4 Y^4+\lambda_5 Y^3 Z) \Big)=0,
\label{3e18}\\
&K_1=
\left|
\begin{array}{cccc}
-\lam_{0} &  \lam_1/2 & 2Z  & -2Y\\
\lam_1/2      & -\lam_2-4Z        & \lam_3/2+2Y & 2X\\
 2Z    & \lam_3/2+2Y &   -\lam_4-4X     & \lam_5/2 \\
-2Y  &  2X   & \lam_5/2   & 0 \\
\end{array}
\right| .
\label{3e19}
\end{align}
In the $\lambda_6=0$ case, the generalized Kummer surface 
relation $K_2=0$ reduce to the standard Kummer surface relation 
$K_1=0$~\cite{Baker2}.

If we consider the special case of $\lambda_6=0$, we obtain
$\displaystyle{\frac{\partial \wp_{21}}{\partial u_1}=\frac{\partial Q}{\partial u_2}}$, 
then we can put $Q=\wp_{11}$. Thus, we obtain conventional differential equations for 
Weierstrass-type hyperelliptic functions~\cite{Baker2}. 
If we consider the more special case of $\lambda_6=0, \lambda_0=0$, 
the differential equations for the Weierstrass-type hyperelliptic differential functions become
\begin{align}
&1)\ \wp_{2222}-6 \wp_{22}^2-\lambda_4 \wp_{22}- \lambda_5 \wp_{21}
-\frac{\lambda_3 \lambda_5}{8} 
 =0,
\label{3e20}\\
&2)\ \wp_{2221}
-6 \wp_{22}\wp_{21}-\lambda_4 \wp_{21}+\frac{\lambda_5}{2} \wp_{11}
=0,
\label{3e21}\\
&3)\ \wp_{2211}
-2 \wp_{22} \wp_{11}-4 \wp_{21}^2
-\frac{\lambda_3}{2} \wp_{21}=0,
\label{3e22}\\
&4)\ \wp_{2111}
-6 \wp_{21} \wp_{11}
+\frac{\lambda_1}{2} \wp_{22}-\lambda_2 \wp_{21}=0,
\label{3e23}\\
&5)\ \wp_{1111}
-6 \wp_{11}^2 -\lambda_1 \wp_{21}-\lambda_2 \wp_{11}
-\frac{\lambda_1 \lambda_3}{8}=0 .
\label{3e24}
\end{align}
The condition $\lambda_0=0$ is not a special  restriction, 
because we can always put $\lambda_0=0$ by the constant shift of $x_i$ in the 
hyperelliptic curve. 
\subsection{Genus two Jacobi type hyperelliptic differential equations} 
Next, we try to obtain differential equations for the genus two Jacobi-type hyperelliptic 
functions using the following dual transformation.
Under a dual transformation 
$$x_i \rightarrow \frac{1}{x_i}, \quad y_i \rightarrow \frac{y_i}{x_i^3}, \quad 
\lambda_j \leftrightarrow \lambda_{6-j},\ (j=0,1, \cdots , 6),$$
hyperelliptic curves of Eq.(\ref{3e1}) are invariant. Under this transformation, we obtain 
\begin{align}
&du_1=\frac{dx_1}{y_1}+\frac{dx_2}{y_2} \rightarrow 
-\left(\frac{x_1 dx_1}{y_1}+\frac{x_2dx_2}{y_2}\right)=-du_2 , 
\label{3e25}\\
&du_2=\frac{x_1 dx_1}{y_1}+\frac{x_2 dx_2}{y_2} \rightarrow
-\left(\frac{dx_1}{y_1}+\frac{dx_2}{y_2}\right)=-du_1 .  
\label{3e26}
\end{align}
Then the corresponding hyperelliptic functions become of the form 
\begin{align}
&\wp_{22}(u_1,u_2)=\frac{\lambda_5}{4} (x_1+x_2)  \rightarrow 
\wwp_{11}(u_1, u_2)=\frac{\lambda_1}{4}\left( \frac{1}{x_1}+\frac{1}{x_2}\right),
\label{3e27}\\
&\wp_{21}(u_1,u_2)=\frac{\lambda_5}{4} x_1x_2  \rightarrow 
\wwp_{21}(u_1, u_2)=-\frac{\lambda_1}{4}\frac{1}{ x_1 x_2},
\label{3e28}\\
& Q(u_1, u_2) \rightarrow 
\ \Q(u_1, u_2)=\frac{F(x_1,x_2;\left\{\lambda_j\right\})-2 y_1 y_2}{4 x_1 x_2 (x_1-x_2)^2}.
\label{3e29}
\end{align}
We call these hyperelliptic functions $\wwp_{11}$, $\wwp_{21}$ and $\Q$ as the Jacobi type.
These satisfy relations 
\begin{eqnarray}
\frac{ \partial \Q}{ \partial u_1}
\ne \frac{ \partial \wwp_{21}}{ \partial u_2}, 
\quad \frac{ \partial \wwp_{21}}{ \partial u_1}
=\frac{ \partial \wwp_{11}}{ \partial u_2} . 
\label{3e30}
\end{eqnarray}
%
Differential equations for the genus two Jacobi type hyperelliptic functions are given by 
\begin{align}
&1)\ \wwp_{2221}-6 \wwp_{21} \Q -\lambda_4 \wwp_{21} 
+\frac{\lambda_5}{2} \wwp_{11} -\frac{32 \lambda_0}{\lambda_1^2} \wwp_{21}^3
+\frac{\lambda_1 \lambda_6}{4}=0 , 
\label{3e31}\\
&2)\ \wwp_{2211}-2\wwp_{11} \Q -4 \wwp_{21}^2-\frac{\lambda_3}{2}\wwp_{21}
-\frac{32 \lambda_0}{\lambda_1^2} \wwp_{21}^2 \wwp_{11} =0,
\label{3e32}\\
&3)\ \wwp_{2111}-6 \wwp_{21} \wwp_{11}+\frac{\lambda_1}{2} \Q
-\lambda_2 \wwp_{21} -\frac{4 \lambda_0}{\lambda_1} \wwp_{21}^2
-\frac{32 \lambda_0}{\lambda_1^2} \wwp_{21} \wwp_{11} ^2=0,
\label{3e33}\\
&4)\ \wwp_{1111}-6 \wwp_{11}^2-\lambda_1 \wwp_{21}-\lambda_2 \wwp_{11}
-\frac{12 \lambda_0}{\lambda_1} \wwp_{21} \wwp_{11}
-\frac{32 \lambda_0}{\lambda_1^2} \wwp_{11} ^3-\frac{\lambda_1 \lambda_3}{8}=0,
\label{3e34}\\
&5)\ \Q_{11}-2 \wwp_{11} \Q 
-\left(4-\frac{16 \lambda_0 \lambda_2}{\lambda_1^2}\right) \wwp_{21}^2
+\frac{\lambda_0 \lambda_5}{\lambda_1} \wwp_{11}
-\left(\frac{\lambda_3}{2}+\frac{2\lambda_0 \lambda_4}{\lambda_1}\right) \wwp_{21}
\nonumber\\
&-\frac{20 \lambda_0}{\lambda_1} \wwp_{21} \Q
-\frac{8 \lambda_0 \lambda_3}{\lambda_1^2} \wwp_{21} \wwp_{11}
-\frac{32 \lambda_0}{\lambda_1^2} \wwp_{11} ^2 \Q
+\frac{\lambda_0 \lambda_6}{2}=0,
\label{3e35}\\
&6)\ \Q_{21}-6 \wwp_{21} \Q -\lambda_4 \wwp_{21} 
+\left(\frac{\lambda_5}{2}+\frac{2\lambda_0 \lambda_6}{\lambda_1}\right) \wwp_{11}
+\frac{4 \lambda_0 \lambda_3}{\lambda_1^2} \wwp_{21}^2
\nonumber\\
&-\frac{8 \lambda_0 \lambda_4}{\lambda_1^2} \wwp_{21} \wwp_{11}
+\frac{4 \lambda_0 \lambda_5}{\lambda_1^2} \wwp_{11}^2
-\frac{32 \lambda_0}{\lambda_1^2} \wwp_{21} \wwp_{11} \Q
+\frac{\lambda_1 \lambda_6}{4}=0,
\label{3e36}\\
&7)\ \Q_{22} -6 \Q^2-\lambda_4 \Q +3 \lambda_6 \wwp_{11}
-\left(\lambda_5-\frac{2\lambda_0 \lambda_6}{\lambda_1}\right) \wwp_{21}
-\frac{8 \lambda_0 \lambda_5}{\lambda_1^2} \wwp_{21} \wwp_{11}
\nonumber\\
&+\frac{16 \lambda_0 \lambda_6}{\lambda_1^2} \wwp_{11}^2
-\frac{32 \lambda_0 }{\lambda_1^2} \wwp_{21}^2 \Q
+\left(\frac{\lambda_2 \lambda_6}{2}-\frac{\lambda_3 \lambda_5}{8}\right)=0 .
\label{3e37}
\end{align}
From Eq.(\ref{3e33}), $\Q$ naturally appears as a combination of 
the terms in Eq.(\ref{3e29}).
The characteristic feature of the $\lambda_0 \ne 0$ hyperelliptic differential equations is the existence
 of the cubic terms  
$\displaystyle{-\frac{32 \lambda_0 }{\lambda_1^2} \wwp_{21}^3,
-\frac{32 \lambda_0 }{\lambda_1^2} \wwp_{21}^2 \wwp_{11}, \cdots, 
-\frac{32 \lambda_0 }{\lambda_1^2} \wwp_{21}^2  \Q}$,  
 in Eq.({\ref{3e31}) $\sim$ Eq.(\ref{3e37}).\\
If we take the special case of $\lambda_6=0$, differential equations  
for the genus two Weierstrass type hyperelliptic functions does not 
contain cubic terms just as the differential equation of the genus one Weierstrass 
$\wp$ function.
The above differential equations for genus two Jacobi-type hyperelliptic functions 
contain cubic terms just as the differential equation of the genus one 
Jacobi $\sn$ function. This is why we call $\wwp_{i j}$, $\Q$ Jacobi-type
hyperelliptic functions.

In the special case of $\lambda_0=0, \lambda_6=0$, we obtain 
$\displaystyle{\frac{ \partial \Q}{ \partial u_1}=\frac{ \partial \wwp_{21}}{ \partial u_2}}$, 
then we can put $\Q=\wwp_{22}$ in such a special case. 
Therefore, the correspondence in Eq.(\ref{3e29}) is as follows: 
$\wp_{22}\rightarrow \wwp_{11}, \wp_{21}\rightarrow \wwp_{21}, \wp_{11}\rightarrow \wwp_{22}$.
In this special case, differential equations for the genus 
two Jacobi type hyperelliptic equations become
\begin{align}
&1)\ \wwp_{2222}-6 \wwp_{22}^2-\lambda_4 \wwp_{22}- \lambda_5 \wwp_{21}
-\frac{\lambda_3 \lambda_5}{8} 
 =0,
\label{3e38}\\
&2)\ \wwp_{2221}
-6 \wwp_{22}\wwp_{21}-\lambda_4 \wwp_{21}+\frac{\lambda_5}{2} \wwp_{11}
=0,
\label{3e39}\\
&3)\ \wwp_{2211}
-2 \wwp_{22} \wwp_{11}-4 \wwp_{21}^2
-\frac{\lambda_3}{2} \wwp_{21}=0,
\label{3e40}\\
&4)\ \wwp_{2111}
-6 \wwp_{21} \wwp_{11}
+\frac{\lambda_1}{2} \wwp_{22}-\lambda_2 \wwp_{21}=0,
\label{3e41}\\
&5)\ \wwp_{1111}
-6 \wwp_{11}^2 -\lambda_1 \wwp_{21}-\lambda_2 \wwp_{11}
-\frac{\lambda_1 \lambda_3}{8}=0.
\label{3e42}
\end{align}
Therefore, in the special case of $\lambda_0=0, \lambda_6=0$, the genus two Weierstrass type 
hyperelliptic functions $\{ \wp_{22}. \wp_{21}, \wp_{11} \}$
and the Jacobi type hyperelliptic functions $\{ \wwp_{22}. \wwp_{21}, \wwp_{11} \}$ satisfy 
same set of differential equations. 
The relation between two solutions is given by
\begin{align}
\wwp_{22}=-\frac{\lambda_5}{4} \frac{ \wp_{11}}{\wp_{21}} , \ 
\wwp_{21}=\frac{\lambda_1 \lambda_5}{16} \frac{1}{\wp_{21}} , \ 
\wwp_{11}=-\frac{\lambda_1}{4} \frac{\wp_{22}}{\wp_{21}} ,
\label{3e43}
\end{align}
which corresponds the half-period relation 
$$(\wp(u+\omega_3)-e_3)=\frac{(e_3-e_1)(e_3-e_2)}{(\wp(u)-e_3)},  \quad
\sn(u+3iK')=\frac{1}{k\ \sn(u)}, $$
in genus one case.
Eq.(\ref{3e43}) is a special case of the half-period relation of 
the standard projective type II ${\rm Sp(4,{\bf R})}$ Lie group
transformation 
\begin{eqnarray}
\wp_{22}(u+\Omega )=\frac{a_1 \wp_{22}(u) 
+a_2 \wp_{21}(u) +a_3 \wp_{11}(u)+a_4}
{d_1 \wp_{22}(u) +d_2 \wp_{21}(u) +d_3 \wp_{11}(u)+d_4}  ,
\label{3e44}\\
\wp_{21}(u+\Omega)=\frac{b_1 \wp_{22}(u) 
+b_2 \wp_{21}(u) +b_3 \wp_{11}(u)+b_4}
{d_1 \wp_{22}(u) +d_2 \wp_{21}(u) +d_3 \wp_{11}(u)+d_4}  ,
\label{3e45}\\
\wp_{11}(u+\Omega)=\frac{c_1 \wp_{22}(u) 
+c_2 \wp_{21}(u) +c_3 \wp_{11}(u)+c_4}
{d_1 \wp_{22}(u) +d_2 \wp_{21}(u) +d_3 \wp_{11}(u)+d_4}, 
\label{3e46}
\end{eqnarray}
and
\begin{equation}
G_{{\rm II}}=\left( \begin{array}{@{\,}cccc@{\,}} 
a_1 & a_2 & a_3 & a_4 \\ 
b_1 & b_2 & b_3 & b_4 \\ 
c_1 & c_2 & c_3 & c_4 \\ 
d_1 & d_2 & d_3 & d_4 \\ 
\end{array}\right)
=\left( \begin{array}{@{\,}cccc@{\,}} 
0 & 0  &  1  &  0\\ 
0  & 0 & 0 &  -1\\ 
1 & 0 & 0 & 0 \\ 
0 & -1 &  0 & 0 \\ 
\end{array}\right)  .
\label{3e47}
\end{equation}
by putting $\lambda_5=4, \lambda_1=4$~\cite{H1}.
 
\section{Summary and Discussions} 
\setcounter{equation}{0}
The static KdV equation has solutions of the Weierstrass $\wp$ function.
The Weierstrass $\wp$ elliptic function and Jacobi $\sn$ elliptic functions 
form a family of elliptic functions.
We have generalized the mKdV equation in such a way that 
the static mKdV differential equation 
and static KdV differential equation have good correspondence.
We demonstrated that this generalized mKdV equation is integrable. 
We have also demonstrated that there are four kinds of relations between the solutions 
of the static KdV equation and the generalized static mKdV equation.
Next, we generalized the relationship between $\wp$ and $\sn$ functions 
to the genus two case.
 
For the general genus two hyperelliptic curves, we obtained differential equations 
for the genus two Weierstrass-type and Jacobi-type hyperelliptic functions.

If we consider the special case of $\lambda_6=0, \lambda_0=0$, the differential equations 
of the genus two Weierstrass-type and Jacobi-type 
hyperelliptic functions become the same. Then the Weierstrass-type and the 
Jacobi-type hyperelliptic functions are linked as half-period shift solutions 
of the same differential equations. In fact, we confirmed that two hyperelliptic 
functions are linked by a special half-period ${\rm Sp(4, {\bf R})}$ 
Lie group transformation.


\end{document}